\documentclass[12pt]{article}
\usepackage{amsmath,amssymb,amsfonts}
\usepackage{color}
\numberwithin{equation}{section}
\newtheorem{thm}{Theorem}[section]

\newtheorem{prop}[thm]{Proposition}
\newtheorem{lem}[thm]{Lemma}

\newtheorem{defn}{Definition}

\newtheorem{que}{Question}
{ \pagestyle{plain}

\def\cwedge{\bigcirc\kern-1.07em\wedge\ }

\newcommand{\qed}{\hfill\fbox{}\par\vspace{.2cm}}

\begin{document}

\begin{center}
{\LARGE \bf  {{A generalization of contact metric manifolds}}}
\end{center}
{
\begin{center}
{\large \bf J. H. Kim$^{1}$, J. H. Park$^{1} $ and K.
Sekigawa$^{2}$}\end{center} }
\begin{center}
$^{1}$Sungkyunkwan University,
 Suwon, Korea\\
$^{2}$Niigata University,
    Niigata, Japan\end{center}

%\title{\bf{{A generalization of contact metric manifolds}} }
%\author{JangHyen Kim$^{*}$, Jeong Hyeong Park$^{*}$, and K. Sekigawa$^{\dag}$
%\thanks{Department of Mathematics, Sungkyunkwan University,
% Suwon 440-746, Korea, e-mail: parkj@skku.edu}, \thanks{Department of Mathematics,
%    Niigata University,
%    Niigata 950-2181, JAPAN, e-mail: sekigawa@math.sc.niigata-u.ac.jp}
 %}
%
%\date{}

% \maketitle

%%%%%%%%%%%%%%%%%%%%%%%%%%%%%%%%%%%%%%%%%%%%%%%%%%%%%%%%%%%%%%%%%%%%%%%%%
\begin{abstract}
We give a characterization of a contact metric manifold as a special
almost contact metric manifold and discuss an almost contact metric
manifold which is {a} natural generalization of the contact metric
manifolds introduced by Y. Tashiro.
\end{abstract}
\noindent {\it Mathematics Subsect Classification (2010)} : 53B20, 53C20

%%%%%%%%%%%%%%%%%%%%%%%%%%%%%%%%%%%%%%%%%%%%%%%%%%%%%%%%%%%%%%%%%%%%%%%%%%%%%%%%%
\section{Introduction}\label{sec1}
%%%%%%%%%%%%%%%%%%%%%%%%%%%%%%%%%%%%%%%%%%%%%%%%%%%%%%%%%%%%%%%%%%%%%%%%%%%%%%%%%

\indent A $(2n+1)$-dimensional smooth manifold $M$ is called a {\it
contact manifold} if it admits a global $1$-form $\eta$ such that {{
$\eta\wedge(d\eta)^{n}\neq 0$ everywhere on $M$. Then we call the 1-form $\eta$ a contact
form of $M$. It is well-known
that given a contact form $\eta$, there exists a unique vector field
$\xi$, which is called the {\it characteristic vector field},
satisfying $\eta(\xi)=1$ and $d\eta(\xi,X)=0$ for any vector field
$X$ on $M$. A Riemannian metric $g$ is said to be an
\emph{associated metric} to a contact form $\eta$ if there exists a
$(1,1)$-tensor field $\phi$ satisfying
\begin{equation}\label{11}
        \eta(X)=g(X,\xi),\quad d\eta(X,Y)=g(X,\phi Y)
%        \quad \phi^2 X=- X+\eta( X)\xi
\end{equation}
for $X, Y  \in\mathfrak{X}(M)$. A $(2n+1)$-dimensional smooth
manifold equipped with a triple $(\phi, \xi, \eta)$ of a
$(1,1)$-tensor field $\phi$, a vector field $\xi$ and a 1-form
$\eta$ on $M$ satisfying

%From \eqref{21}, one can easily obtain
\begin{equation}\label{12}
         \phi^2 X=- X+\eta( X),\quad \phi\xi=0,\quad
         \eta\circ\phi=0,\ \ and \ \ \eta(\xi)=1
%        \quad g(\phi X,\phi Y)=g(X,Y)-\eta(X)\eta(Y).
\end{equation}
for $X  \in\mathfrak{X}(M)$ is called an {\it {almost contact
manifold}} with the almost contact structure $(\phi, \xi, \eta)$.
Further, an almost contact manifold $M=(M, \phi, \xi, \eta)$
equipped with a Riemannian metric $g$ satisfying

\begin{equation}\label{13}
      g(\phi X,\phi Y)=g(X,Y)-\eta(X)\eta(Y),\quad \eta(X)=g(\xi,X)
\end{equation}
for $X, Y  \in\mathfrak{X}(M)$ is called an {\it{almost contact
metric manifold}} with the almost contact metric structure $(\phi,
\xi, \eta, g)$. From (\ref{11}) $\sim$ (\ref{13}), we may regard a
contact metric manifold as a special contact metric manifold.

% The structure $(\phi, \xi, \eta, g)$ is called a {\it contact metric
% structure}, and a manifold $M$ with a contact metric structure
% $(\phi, \xi, \eta, g)$ is said to be a \emph{contact metric
% manifold} and is denoted by $M = (M, \phi, \xi,\eta, g)$.

% It is well-known that a $(2n+1)$-dimensional smooth manifold $M$ admits an
% almost contact structure if and only if the structure group of the
% tangent bundle $TM$ of its reducible to the group $U(n)\times 1$.
% This fact is corresponding tn the one that a $2n$-dimensional smooth
% manifold $\bar{M}$ admits an almost complex structure if and only if
% the structure group the tangent bundle $T\bar{M}$ of $\bar{M}$ is
% reducible to $U(m)$.

 D. Chinea and C. Gonzalez \cite{CG} obtained a classification of
the $(2n+1)$-dimensional almost contact metric manifold based on
$U(n) \times I$ representation theory, which is an analogy of the
classification of the $2n$-dimensional almost Hermitian manifolds
established by A. Gray and H. M. Hervella \cite{GH}.

Now, let $M=(M, \phi , \xi, \eta , g)$ be a $(2n+1)$-dimensional
almost contact metric manifold and $\bar{M}=M \times {\mathbb{R}}$
be the product manifold of  $M$ and a real line ${\mathbb{R}}$
equipped with the following almost Hermitian structure
$(\bar{J},\bar{g})$ defined by

\begin{equation}\label{14}
\begin{split}
& \bar{J}X=\phi X-\eta(X)\frac{\partial}{\partial{t}},
   \quad \bar{J}\frac{\partial}{\partial{t}}=\xi,\\
& \bar{g}(X,Y)=e^{-2t}g(X,Y),
   \quad \bar{g}(X,\frac{\partial}{\partial{t}})=0,\\
& \bar{g}(\frac{\partial}{\partial{t}},\frac{\partial}{\partial{t}})=e^{-2t}\\
\end{split}
\end{equation}
for $X, Y  \in\mathfrak{X}(M)$ and $t \in {\mathbb{R}}$. In the case
where $\bar{J}$ is integrable, the corresponding almost contact
metric manifold $M=(M, \phi, \xi, \eta, g)$ is said to be
\emph{normal}. Especially, a normal contact metric manifold is
called a Sasaki manifold. Y. Tashiro \cite{T2} discussed the
relation ship between the classes of almost Hermitian manifolds and
the corresponding ones of
almost contact metric manifolds and showed the following:\\

% structures and {\red the one }of the corresponding almost
% Hermitian manifolds and {\red proved }the followings :

Fact 1. $\bar{M}=(\bar{M},\bar{J},\bar{g})$ is a K\"ahler manifold
if and only if $M=(M, \phi, \xi, \eta, g )$  is a Sasakian manifold.

Fact 2. $\bar{M}=(\bar{M},\bar{J},\bar{g})$  is an almost K\"ahler
manifold if  and only  if $M=(M, \phi, \xi, \eta, g
)$ is a contact  metric manifold.\\

On the other hand, it is easily observed that any orientable
hypersurface of an almost Hermitian manifold becomes an almost
contact metric manifold in natural way. So, from the above
observation, it seems natural to consider the almost contact metric
manifold in connection with almost Hermitian geometry, for example,
to
discuss the classification  %following the Gray-Hervella
of almost Hermitian manifolds. We denote by $\mathcal{K}$,
$\mathcal{AH}$, $\mathcal{NK}$, $\mathcal{QK}$ and $\mathcal{H}$ the
classes of K\"ahler manifolds, almost K\"ahler manifolds, nearly
K\"ahler manifolds, quasi K\"ahler manifolds and Hermitian
manifolds, respectively {thus, their inclusion} relations are as
follows \cite{GH}:
%%%%
% {\blue We denote by $\mathcal{AH}$, $\mathcal{H}$, $\mathcal{K}$,
% $\mathcal{AK}$, $\mathcal{NK}$, and $\mathcal{QK}$ the sets of all
% almost Hermitian manifolds, Hermitian manifolds, K\"ahler manifolds,
% almost K\"ahler manifolds, nearly K\"ahler manifolds, and quasi
% K\"ahler manifolds,
% respectively. Then the inclusion relations among them are as
% follows :
%%%%
\begin{equation}\label{15}
%$$
    \mathcal K
    \begin{matrix}
    \subset \mathcal{AK} \subset\\
    \subset \mathcal{NK} \subset
    \end{matrix}
    \mathcal{QK},
    % \subset \mathcal{SK},
   \qquad
    %\mathcal K = \mathcal{QK} \mathcal {H},
   % \qquad
   % \mathcal K = \mathcal{NK} \cap mathcal{AK}.
\end{equation}
%$\mathcal{K}\subset \mathcal{QK}$
\begin{equation*}
\begin{split}
&\mathcal{AK}\cap \mathcal{NK} = \mathcal{K},\\
&\mathcal{QK}\cap \mathcal{H} = \mathcal{K}.
\end{split}
\end{equation*}
}

%%%
%\begin{equation}\label{14}
%\begin{split}
%K & \subset AK \subset QK,\\
%  & \subset NK \\
%AK & \cap NK = K ,\quad QK \cap H=H.
%\end{split}
%\end{equation}

\noindent In the next section, we shall reprove these facts and
introduce a class of almost contact metric manifolds as the class of
almost contact metric manifolds corresponding to the class of quasi
K\"ahler manifolds, which is regarded as a generalization of the
class of contact metric manifolds by taking account of (\ref{15}).

In the sequel, we shall call such an almost contact metric manifold
quasi contact metric manifold. In $\S4$, we shall discuss the quasi
contact metric manifolds from the view point of a generalization of
contact metric manifolds.

% In the nest section, we shall reprove these facts. \ From the above
% observation, {\red it seems ako ororthmhule the classify} all almost
% contact metric manifolds {\red accordingly} the classification
% $M=(M, \phi, \xi, \eta, g )$ of the corresponding almost Hermitian
% manifolds $\bar{M}=(\bar{M},\bar{J},\bar{g})$.

% {\red Taking account of the Gray-Hervella classification theorem},
% it is easily checked that a quasi K$\ddot{a}$hler manifold a
%% generaligation of both of almost K$\ddot{a}$hler manifold and a
% nearly K$\ddot{a}$hler manifold.

%%%%%%%%%%%%%%%%%%%%%%%%%%%%%%%%%%%%%%%%%%%%%%%%%%%%%%%%%%%%%%%%%%%%%%%%%%%%%%%%%
\section{Preliminaries}\label{sec2}
%%%%%%%%%%%%%%%%%%%%%%%%%%%%%%%%%%%%%%%%%%%%%%%%%%%%%%%%%%%%%%%%%%%%%%%%%%%%%%%%%

In this section, we shall prepare some fundamental formulas which we
need in the forthcoming discussions in the present paper. Let $M=(M,
\phi, \xi, \eta, g)$ be a $(2n+1)$-dimensional almost contact metric
manifold and $\bar{M}=M \times {\mathbb{R}}$  be the direct product
manifold of $M$ and a real line equipped with the almost Hermitian
structure $(\bar{J},\bar{g})$ defined by (1.4). Now, we denote by
$[\phi,\phi]$ the (1,2)-tensor field defined by

\begin{equation}\label{21}
[\phi,\phi](X,Y)=[\phi X,\phi Y]-[X,Y]-\phi[\phi X, Y]-\phi[X,\phi
Y]+\eta([X,Y])\xi
\end{equation}
for $X, Y  \in\mathfrak{X}(M)$. Further, we denote by $\bar{N}$ the
Nijenhuis tensor of the almost complex structure $\bar{J}$. Then,
from (\ref{14}), we have

\begin{equation}\label{22}
\bar{N}(X,Y)=[\phi,\phi](X,Y)+2d\eta(X,Y)\xi -\bigg( (L_{\phi X}
\eta)(Y)-(L_{\phi Y} \eta)(X) \bigg) \frac{\partial}{\partial t},
\end{equation}

\begin{equation}\label{23}
\bar{N}(X,\frac{\partial}{\partial t})=-(L_{\xi} \phi)X+(L_{\xi}
\eta)(X)\frac{\partial}{\partial t}
\end{equation}
for $X, Y  \in\mathfrak{X}(M)$. We denote by $N^{(1)}$, $N^{(2)}$,
$N^{(3)}$ and $N^{(4)}$ the following tensor fields on $M$ defined
respectively by

\begin{equation}\label{24}
N^{(1)}(X,Y)=[\phi,\phi](X,Y)+2d\eta(X,Y)\xi,
\end{equation}

\begin{equation}\label{25}
N^{(2)}(X,Y)=(L_{\phi X} \eta)(Y)-(L_{\phi Y} \eta)(X),
\end{equation}

\begin{equation}\label{26}
N^{(3)}(X)=-(L_{\xi} \phi)X,
\end{equation}

\begin{equation}\label{27}
 N^{(4)}(X)=(L_{\xi} \eta)(X)
\end{equation}
for $X, Y  \in\mathfrak{X}(M)$. then, from (\ref{22}) $\sim$
(\ref{27}), we have

\begin{equation}\label{28}
\begin{split}
&\bar{N}(X,Y)=N^{(1)}(X,Y){ -
}N^{(2)}(X,Y)\frac{\partial}{\partial t},\\
%\end{equation}
%\begin{equation}\label{29}
&\bar{N}(X,\frac{\partial}{\partial
t})=N^{(3)}(X)+N^{(4)}(X)\frac{\partial}{\partial t}
\end{split}
\end{equation}
for $X, Y  \in\mathfrak{X}(M)$.

\begin{prop}\label{2.1}\cite{B1}
For an almost contact manifold $M=(M,\phi, \xi, \eta)$ the vanishing
of the tensor field $N^{(1)}$ implies the vanishing of the tensor
fields $N^{(2)}$, $N^{(3)}$ and $N^{(4)}$.
\end{prop}

\begin{prop}\label{2.2}\cite{B1}
For a contact metric manifold $M=(M,\phi,\xi,\eta,g)$, $N^{(2)}$ and
$N^{(4)}$ vanish. Moreover, $N^{(3)}$ vanishes if and only if $\xi$
is a Killing vector field (namely, $M$ is a K-contact manifold).
\end{prop}

\noindent {\bf{Remark}} \ From Proposition \ref{2.1}, taking account
of (\ref {28}), we see that an almost contact metric manifold
$M=(M,\phi,\xi,\eta,g)$ is normal if and only if $N^{(1)}$ vanishes
everywhere on $M$ [1, p.71].

\noindent We here note that the following equality
\begin{equation}\label{210}
\begin{split}
N^{(2)}(X,Y) &= (L_{\phi X}\eta)(Y)-(L_{\phi Y}\eta)(X)\\
             &= \phi X(\eta(Y))-\eta([\phi X,Y]) - \phi Y(\eta(X))+
                \eta([\phi Y,X])\\
             &= (\nabla_{\phi X}\eta)(Y)+\eta(\nabla_{\phi X}Y)-
               \eta(\nabla_{\phi X}Y- \nabla_{Y}(\phi X))\\
             & \quad  - (\nabla_{\phi Y}\eta)(X)-\eta(\nabla_{\phi Y}
             X)+\eta(\nabla_{\phi Y} X - \nabla_{X}(\phi Y))\\
             &= (\nabla_{\phi X}\eta)(Y) + \eta( \nabla_{Y}(\phi
               X))-(\nabla_{\phi Y} \eta)(X)-\eta(\nabla_{X}(\phi
               Y))\\
             &= (\nabla_{\phi X}\eta)(Y)-(\nabla_{Y}\eta)(\phi
             X)-(\nabla_{\phi Y}\eta)(X)+(\nabla_{X}\eta)(\phi Y)\\
\end{split}
\end{equation}
for $X, Y  \in\mathfrak{X}(M)$. We here define a (1,1)-tensor field
$h$ on $M$ by

\begin{equation}\label{211}
h=\frac{1}{2}L_{\xi} \phi.
\end{equation}

\noindent The tensor field $h$ plays an important role in the
geometry of almost contact metric manifolds. From (\ref{211}), we
have easily the following equalities

\begin{equation}\label{212}
hX=\frac{1}{2} \bigg((\nabla_{\xi} \phi)X - \nabla_{\phi X}\xi +
\phi \nabla_{X}\xi \bigg),
\end{equation}

\noindent and hence
\begin{equation}\label{213}
h \xi=0,
\end{equation}

\begin{equation}\label{214}
tr h=0.
\end{equation}

\begin{prop}\label{2.3}
Let $M=(M,\phi,\xi,\eta,g)$ be an almost contact metric manifold
satisfying $\nabla_{\xi} \phi=0$. Then $h$ is symmetric with respect
to the metric $g$ if and only if $N^{(2)}$ vanishes everywhere on
$M$.
\end{prop}

\noindent{\bf{Proof}}. By the hypothesis from (\ref{210}) and
(\ref{212}), we have

\begin{equation}\label{215}
\begin{split}
g(hX,Y)-g(X,hY) & = \frac{1}{2}\bigg( -(\nabla_{\phi
                      X}\eta)(Y)-(\nabla_{X}\eta)(\phi Y) + (\nabla_{\phi Y}\eta)(X) +
                        (\nabla_{Y}\eta)(\phi X) \bigg)\\
               & = - \frac{1}{2}N^{(2)}(X,Y)
\end{split}
\end{equation}
for $X, Y  \in\mathfrak{X}(M)$. This Proposition \ref{2.3} follows
immediatily from (\ref{215}). The following is well-known. \; \; \;
\;\;\;\;\qquad\qquad\qquad\qquad\qquad\qquad\qquad\qquad\qquad\qquad\qquad\qquad\qquad\;\;
$\square$

\begin{prop}\label{2.4}
An almost contact metric manifold $M=(M,\phi,\xi , \eta, \xi, \eta)$
is Sasakian if and only if
$(\nabla_{X}\phi)Y=g(X,Y)\xi-\eta(Y)X$ holds for any $ X, Y
\in\mathfrak{X}(M)$.
\end{prop}

\noindent Now, we denote by $\bar{\nabla}$ the covariant derivative
with respect to the metric $\bar{g}$ on $\bar{M}$. Then, from (1.4)
by direct calculation, we have

\begin{equation}\label{216}
\begin{split}
& \bar{\nabla}_{X}Y =\nabla_{X}Y+g(X,Y)\frac{\partial}{\partial t},\\
& \bar{\nabla}_{X}\frac{\partial}{\partial t} = -X,\\
& \bar{\nabla}_{\frac{\partial}{\partial t}}X = -X,\\
& \bar{\nabla}_{\frac{\partial}{\partial t}}\frac{\partial}{\partial
t} = - \frac{\partial}{\partial t}
\end{split}
\end{equation}
for $X, Y  \in\mathfrak{X}(M)$. Thus, from (\ref{14}) and
(\ref{216}), we have further

\begin{equation}\label{217}
(\bar{\nabla}_{X}\bar{J})Y=(\nabla_{X}\phi)Y-g(X,Y)\xi+\eta(Y)X -
\bigg( g(\phi X,Y)+(\nabla_{X}\eta)(Y)
\bigg)\frac{\partial}{\partial t},
\end{equation}

\begin{equation}\label{218}
(\bar{\nabla}_{X}\bar{J})\frac{\partial}{\partial t}=
\nabla_{X}\xi+\phi X,
\end{equation}

\begin{equation}\label{219}
(\bar{\nabla}_{\frac{\partial}{\partial t}}\bar{J})X=0, \quad
(\bar{\nabla}_{\frac{\partial}{\partial
t}}\bar{J})\frac{\partial}{\partial t}=0
\end{equation}
for $X, Y  \in\mathfrak{X}(M)$. We here show the Facts 1 and 2, from
(\ref{217}) $\sim$ (\ref{219}), we see that
$\bar{M}=(\bar{M},\bar{J},\bar{g})$ is K\"ahler if and only if
\begin{equation}\label{220}
(\nabla_{X}\phi)Y-g(X,Y)\xi+\eta(Y)X=0,
\end{equation}
and
\begin{equation}\label{221}
\nabla_{X}\xi+\phi X=0
\end{equation}
for $X, Y  \in\mathfrak{X}(M)$. Thus, from (\ref{220}) and
(\ref{221}), taking account of Proposition \ref{2.4}, it follows
immediately that $M=(M,\phi,\xi,\eta,g)$ is Sasakian.

Similarly, from (\ref{217}) $\sim$ (\ref{219}), taking account of
(\ref{14}), we see that $\bar{M}=(\bar{M},\bar{J},\bar{g})$ is
almost K\"ahler if and only if

\begin{equation}\label{222}
\begin{split}
0 &=\bar{g}((\bar{\nabla}_{X}\bar{J})Y,Z)
    +\bar{g}((\bar{\nabla}_{Y}\bar{J})Z,X)
    +\bar{g}((\bar{\nabla}_{Z}\bar{J})X,Y)\\
  &=e^{-2t} \bigg( g( (\nabla_{X}\phi)Y,Z)
    +g((\nabla_{Y}\phi)Z,X)
    +g( (\nabla_{Z}\phi)X,Y) \bigg)\\
  &=-3e^{-2t}d\Phi(X,Y,Z),
\end{split}
\end{equation}

\begin{equation}\label{223}
\begin{split}
0 &=\bar{g}((\bar{\nabla}_{X}\bar{J})\frac{\partial}{\partial t},Z)
    +\bar{g}((\bar{\nabla}_{\frac{\partial}{\partial t}}\bar{J})Z,X)
    +\bar{g}((\bar{\nabla}_{Z}\bar{J})X,\frac{\partial}{\partial t})\\
  &=e^{-2t} \bigg((\nabla_{X}\eta)(Z)-(\nabla_{Z}\eta)(X)-2\Phi(X,Z) \bigg)\\
\end{split}
\end{equation}
for $X, Y  \in\mathfrak{X}(M)$, where $\Phi(X,Y)=g(X,\phi Y)$. Thus,
from (\ref{222}) and (\ref{223}), it follows that
\begin{equation}\label{224}
\begin{split}
& d\eta(X,Y)=\Phi(X,Y)
\end{split}
\end{equation}
for $X, Y  \in\mathfrak{X}(M)$, and hence $d\Phi=0$. Therefore, we
see that $M=(M,\phi,\xi,\eta,g)$ is a contact metric manifold.

\begin{defn}\label{d1}
An almost contact metric manifold $M=(M,\phi,\xi,\eta,g)$ is called
a quasi contact metric manifold if the corresponding almost
Hermitian manifold $\bar{M}=(\bar{M},\bar{J},\bar{g})$ defined by
(\ref{14}) is a quasi K\"ahler manifold.
\end{defn}

\noindent We note that a quasi contact metric manifold was primary
introduced as a contact $O^*$-manifold by Tashiro \cite{T2}.

Now, we shall derive the condition for an almost contact metric
manifold to be a quasi contact metric manifold. Again, from
(\ref{217}) $\sim$ (\ref{219}), we see that
$\bar{M}=(\bar{M},\bar{J},\bar{g})$  is quasi K\"ahler if and only
if

\begin{equation}\label{225}
\begin{split}
0 = &(\bar{\nabla}_{X}\bar{J})Y+(\bar{\nabla}_{\bar{J}X}\bar{J})\bar{J}Y\\
%  = &(\nabla_{X}\phi)Y-g(X,Y)\xi+\eta(Y)X \\
%    & +(\nabla_{\phi X}\phi)\phi Y -g(X,Y)\xi-\eta(Y)\nabla_{\phi
%    X}\xi + \eta(Y)X \\
%    & - \bigg( (\nabla_{X}\eta)(Y)+(\nabla_{\phi X}\eta)(\phi
%    Y)+2g(\phi X,Y) \bigg)\frac{\partial}{\partial t},\\
  = & (\nabla_{X}\phi)Y + (\nabla_{\phi X}\phi)\phi Y - 2g(X,Y)\xi \\
    & + 2\eta(Y)X - \eta(Y)\nabla_{\phi X}\xi\\
    & - \bigg( (\nabla_{X}\eta)(Y) + (\nabla_{\phi X}\eta)(\phi Y) +
    2g(\phi X,Y) \bigg)\frac{\partial}{\partial t},\\
\end{split}
\end{equation}

\begin{equation}\label{226}
%\begin{split}
0 =  (\bar{\nabla}_{X}\bar{J})\frac{\partial}{\partial t} +
(\bar{\nabla}_{\bar{J}X}\bar{J})\bar{J}\frac{\partial}{\partial t}
  = \nabla_{X}\xi - \phi \nabla_{\phi X}\xi + 2 \phi X,
%\end{split}
\end{equation}

\begin{equation}\label{227}
%\begin{split}
0 =  (\bar{\nabla}_{\frac{\partial}{\partial t}}\bar{J})Y
      +(\bar{\nabla}_{J\frac{\partial}{\partial t}}\bar{J})\bar{J}Y
  =  (\nabla_{\xi}\phi)(\phi Y) - \eta(Y)\nabla_{\xi}\xi - (\nabla_{\xi}\eta)(\phi Y)\frac{\partial}{\partial t}
%\end{split}
\end{equation}
for $X, Y  \in\mathfrak{X}(M)$. Thus, from (\ref{225}) $\sim$
(\ref{227}) it follows that $M=(M,\phi,\xi,\eta,g)$ is a quasi
contact metric manifold if and only if the following equalities

\begin{equation}\label{228}
(\nabla_{X}\phi)Y + (\nabla_{\phi X}\phi)\phi Y = 2g(X,Y)\xi
-2\eta(Y)X+\eta(Y)\nabla_{\phi X}\xi,
\end{equation}

\begin{equation}\label{229}
(\nabla_{X}\eta)(Y) + (\nabla_{\phi X}\eta)(\phi Y) + 2g(\phi
X,Y)=0,
\end{equation}

\begin{equation}\label{230}
(\nabla_{\phi X}\phi)\xi + 2\phi X + \nabla_{X}\xi =0,
\end{equation}

\begin{equation}\label{231}
(\nabla_{\xi}\phi)\phi Y - \eta(Y)\nabla_{\xi}\xi = 0,
\end{equation}

\begin{equation}\label{232}
(\nabla_{\xi}\eta)(\phi Y) = 0.
\end{equation}

\noindent From (\ref{231}), we get set $Y=\xi$, then we get

\begin{equation}\label{233}
\nabla_{\xi}\xi = 0
\end{equation}

\noindent hold for any $ X, Y \in\mathfrak{X}(M)$}. Thus, from
(\ref{231}) and (\ref{233}), we get
%\begin{equation*}
%{\phi}^{2} \nabla_{\xi}\xi = 0,
%\end{equation*}
%and hence

%\begin{equation}\label{234}
%0 = -\nabla_{\xi}\xi + \eta(\nabla_{\xi}\xi)\xi=  -\nabla_{\xi}\xi.
%\end{equation}
\begin{equation}\label{235}
(\nabla_{\xi}\phi)\phi Y = 0.
\end{equation}
From (\ref{235}), we get farther
\begin{equation*}
( \nabla_{\xi}\phi){\phi}^{2}Y = 0,
\end{equation*}
and hence
\begin{equation*}
- ( \nabla_{\xi}\phi)Y + \eta(Y)(\nabla_{\xi}\phi)\xi = 0,
\end{equation*}
and hence
\begin{equation}\label{236}
\nabla_{\xi}\phi = 0.
\end{equation}

\noindent Further, from (\ref{230}), we have
\begin{equation*}
- \phi \nabla_{\phi X}\xi + 2\phi X + \nabla_{X}\xi = 0,
\end{equation*}
and hence
\begin{equation}\label{237}
(\nabla_{\phi X}\eta)(\phi Y)+2g(\phi X,Y) + (\nabla_{X}\eta)(Y) = 0
\end{equation}
for any $X, Y  \in\mathfrak{X}(M)$, which is nothing but
(\ref{229}). Nearly the equality (\ref{230}) is equivalent to the
equality (\ref{229}). Summing up the above arguments, we have the
following:

\begin{prop}\label{2.5}
An almost contact metric manifold $M=(M,\phi,\xi,\eta, g)$ is a
quasi contact metric manifold if and only if the equalities
(\ref{228}), (\ref{229}), (\ref{233}) and (\ref{236}) hold
everywhere on $M$.
\end{prop}

\begin{prop}\label{2.6}\cite{B1}
Let $M=(M,\phi,\xi,\eta, g)$ be an almost contact metric manifold
satisfying the following condition :
\begin{equation*}
(C_{1}) \quad (\nabla_{X}\phi)Y + (\nabla_{\phi X}\phi)\phi Y =
2g(X,Y)\xi - \eta(Y)X - \eta(X)\eta(Y)\xi - \eta(Y)hX
\end{equation*}
for any $ X, Y  \in\mathfrak{X}(M)$. Then, the following equalities
$(C_{2}) \sim (C_{4})$ are derived from the equality $(C_{1})$ :
\begin{equation*}
\begin{split}
& (C_{2}) \quad (\nabla_{X}\eta)Y + (\nabla_{\phi X}\eta)(\phi Y)
+2g(\phi X,Y)=0,\\
%2g(X,Y)\xi - \eta(Y)X - \eta(X)\eta(Y)\xi - \eta(Y)hX,\\
& (C_{3}) \quad \nabla_{\xi}\phi = 0, \\
& (C_{4}) \quad \nabla_{\xi}\xi = 0
\end{split}
\end{equation*}
for any $X, Y  \in\mathfrak{X}(M)$.
\end{prop}

\noindent{\bf{Proof}}. We change $X$ and $Y$ for $\phi X$ and $\phi Y$ in
$(C_{1})$, respectively, we get
\begin{equation*}
(\nabla_{\phi X}\phi)\phi Y +
(\nabla_{(-X+\eta(X)\xi)}\phi)(-Y+\eta(Y)\xi)\\
= 2g(X,Y)\xi -2 \eta(X)\eta(Y)\xi,
\end{equation*}
and hence
\begin{equation*}
\begin{split}
& (\nabla_{ X}\phi)Y + (\nabla_{\phi X}\phi)\phi Y -
\eta(Y)(\nabla_{X}\phi)\xi - \eta(X)(\nabla_{\xi}\phi)Y +
\eta(X)\eta(Y)(\nabla_{\xi}\phi)\xi \\
& = 2g(X,Y)\xi -2 \eta(X)\eta(Y)\xi,
\end{split}
\end{equation*}
namely
\begin{equation}\label{238}
\begin{split}
& (\nabla_{X}\phi)Y + (\nabla_{\phi X}\phi)\phi Y \\
 & =  2g(X,Y)\xi - 2\eta(X)\eta(Y)\xi + \eta(Y)(\nabla_{X}\phi)\xi \\
  & \quad +  \eta(X)(\nabla_{\xi}\phi)Y -\eta(X)\eta(Y)(\nabla_{\xi}\phi)\xi
\end{split}
\end{equation}
for any $X, Y  \in\mathfrak{X}(M)$. Thus, from $(C_{1})$ and
(\ref{238}), we have
\begin{equation*}
\begin{split}
&-\eta(Y)X-\eta(X)\eta(Y)\xi - \eta(Y)hX\\
&= -2\eta(X)\eta(Y)\xi + \eta(Y)(\nabla_{X}\phi)\xi
  +\eta(X)(\nabla_{\xi}\phi)Y-\eta(X)\eta(Y)(\nabla_{\xi}\phi)\xi,
\end{split}
\end{equation*}
and hence
\begin{equation}\label{239}
\begin{split}
& \eta(Y)(\nabla_{X}\phi)\xi +
\eta(X)(\nabla_{\xi}\phi)Y-\eta(X)\eta(Y)(\nabla_{\xi}\phi)\xi \\
& - \eta(X)\eta(Y)\xi + \eta(Y)X + \eta(Y)hX = 0
\end{split}
\end{equation}
for any $X, Y  \in\mathfrak{X}(M)$. Thus, setting $X=Y=\xi$ in
(\ref{239}) and taking account of (\ref{213}), we have
\begin{equation}\label{240}
(\nabla_{\xi}\phi)\xi=0.
\end{equation}
Further, setting $X=\xi$ and choosing $Y$ perpendicular to $\xi$
arbitrary in (\ref{239}), and taking account of (\ref{240}), we have
\begin{equation}\label{241}
(\nabla_{\xi}\phi)Y=0.
\end{equation}
Thus, from (\ref{240}) and (\ref{241}), we have $(C_{3})$.  The
equality $(C_{4})$ follows immediately from $(C_{3})$. Thus, from
(\ref{211}) and $(C_{3})$, we have
\begin{equation}\label{242}
hX=\frac{1}{2}(- \nabla_{\phi X}\xi + \phi \nabla_{X}\xi)
\end{equation}
for $X \in\mathfrak{X}(M)$. Thus from (\ref{239}), taking account of
$(C_{3})$ and (\ref{242}), we obtain
\begin{equation*}
-\eta(Y)\phi \nabla_{X}\xi - \eta(X)\eta(Y)\xi +
\eta(Y)X+\frac{1}{2}\eta(Y)(-\nabla_{\phi X}\xi + \phi
\nabla_{X}\xi) = 0,
\end{equation*}
and hence
\begin{equation}\label{243}
\begin{split}
0 & = \eta(Y) \bigg( -\phi \nabla_{X}\xi + X -\eta(X)\xi
     \frac{1}{2}(-\nabla_{\phi X}\xi + \phi  \nabla_{X}\xi) \bigg) \\
  & = \eta(Y) \bigg( -\frac{1}{2}(\nabla_{\phi X}\xi + \phi
  \nabla_{X}\xi) + X - \eta(X)\xi \bigg)
\end{split}
\end{equation}
for any $X, Y  \in\mathfrak{X}(M)$. Thus, from (\ref{243}), we have
\begin{equation}\label{244}
\nabla_{\phi X}\xi + \phi \nabla_{X}\xi = 2X-2\eta(X)\xi
\end{equation}
for $X \in\mathfrak{X}(M)$. From (\ref{244}), we have also
\begin{equation*}
\phi \nabla_{\phi X}\xi + {\phi}^{2}\nabla_{X}\xi = 2 \phi X,
\end{equation*}
and hence
\begin{equation}\label{245}
-\nabla_{X}\xi + \phi \nabla_{\phi X}\xi = 2 \phi X .
\end{equation}
From (\ref{245}), we have further
\begin{equation*}
(\nabla_{X}\eta)(Y) + (\nabla_{\phi X}\eta)(\phi Y) + 2g(\phi X,Y) =
0
\end{equation*}
for any $X, Y  \in\mathfrak{X}(M)$. Namely, we have $(C_{2})$. \; \;
\; \;\;\;\;\qquad\qquad\qquad\qquad\qquad\qquad\qquad $\square$

In the next section $\S3$, we shall give a characterization for an
almost contact metric manifold to be a contact metric manifold, and
further, a characterization for an almost contact metric manifold to
be a quasi contact metric manifold. Through similar arguments as the
proof of Proposition \ref{2.6}, we have the following:

\begin{prop}\label{2.7}
Let $M=(M, \phi, \xi,\eta, g)$ be an almost contact metric manifold
satisfying the following condition:
\begin{equation*}
(C_{1}') \quad (\nabla_{X}\phi)Y + (\nabla_{\phi X}\phi)\phi Y =
2g(X,Y)\xi -2\eta(Y)X+\eta(Y)\nabla_{\phi X}\xi
\end{equation*}
for any $X, Y  \in\mathfrak{X}(M)$. Then, the equalities
$(C_{2})\sim (C_{4})$ in Proposition \ref{2.6} are derived from
$(C_{1}')$.
\end{prop}

\noindent { \bf{Proof.} } Let $M=(M, \phi, \xi,\eta, g)$ be an
almost contact metric manifold satisfying the condition $(C_{1}')$.
By changing $X$ and $Y$ for $\phi X$ and $\phi Y$ in $(C_{1}')$,
respectively, we get
\begin{equation*}
(\nabla_{\phi X}\phi)\phi Y +
(\nabla_{-X+\eta(X)\xi}\phi)(-Y+\eta(Y)\xi)=2g(X,Y)\xi -
\eta(X)\eta(Y)\xi,
\end{equation*}
and hence
\begin{equation}\label{246}
\begin{split}
 & (\nabla_{X}\phi)Y + (\nabla_{\phi X}\phi)\phi Y \\
          & = 2g(X,Y)\xi -2\eta(X)\eta(Y)\xi +
          \eta(X)(\nabla_{\xi}\phi)Y +
          \eta(Y)(\nabla_{X}\phi)\xi-\eta(X)\eta(Y)(\nabla_{\xi}\phi)\xi
\end{split}
\end{equation}
for any $X, Y  \in\mathfrak{X}(M)$. Thus,$(C_{1}')$ and (\ref{246}),
we have
\begin{equation}\label{247}
\begin{split}
  & -2\eta(X)\eta(Y)\xi +
\eta(Y)(\nabla_{X}\phi)\xi+\eta(X)(\nabla_{\xi}\phi)Y-\eta(X)\eta(Y)(\nabla_{\xi}\phi)\xi
\\
          &= -2\eta(Y)X + \eta(Y)\nabla_{\phi X}\xi
\end{split}
\end{equation}
for any $X, Y  \in\mathfrak{X}(M)$. Setting $X=Y=\xi$ in
(\ref{247}), we have
\begin{equation}\label{248}
(\nabla_{\xi}\phi)\xi = 0.
\end{equation}
Thus, setting $X=\xi,\ Y\perp \xi$ in (\ref{247}), taking account of
(\ref{248}), we have
\begin{equation}\label{249}
(\nabla_{\xi}\phi)Y = 0.
\end{equation}
Thus, from (\ref{248}) and (\ref{249}), we have $(C_{3})$, thus, we
see that (\ref{246}) reduces to
\begin{equation}\label{250}
\begin{split}
 & (\nabla_{X}\phi)Y+(\nabla_{\phi X}\phi)\phi Y \\
%& = 2g(X,Y)\xi-2\eta(X)\eta(Y)\xi  + \eta(Y)(\nabla_{X}\phi)\xi\\
& = 2g(X,Y)\xi-2\eta(X)\eta(Y)\xi  - \eta(Y)\phi \nabla_{X}\xi.
\end{split}
\end{equation}
Thus, from $(C'_{1})$ and (\ref{250}), we have
\begin{equation*}
-2 \eta(Y)X + \eta(Y)\nabla_{\phi X}\xi = -2\eta(X)\eta(Y)\xi -
\eta(Y)\phi \nabla_{X}\xi,
\end{equation*}
and hence
\begin{equation}\label{251}
\eta(Y) \bigg( \nabla_{\phi X}\xi + \phi\nabla_{X}\xi + 2\eta(X)\xi
-2X \bigg) = 0
\end{equation}
for any $X, Y  \in\mathfrak{X}(M)$. From (\ref{251}), we have
further
\begin{equation}\label{252}
\nabla_{\phi X}\xi + \phi\nabla_{X}\xi + 2\eta(X)\xi -2X
    = 0.
\end{equation}
Thus, from (\ref{252}), we have also
\begin{equation*}
\phi \nabla_{\phi X}\xi + {\phi}^{2}\nabla_{X}\xi - 2\phi X = 0,
\end{equation*}
and hence
\begin{equation*}
\phi \nabla_{\phi X}\xi - \nabla_{X}\xi + \eta(\nabla_{X}\xi)\xi
-2\phi X = 0,
\end{equation*}
namely
\begin{equation}\label{253}
\phi \nabla_{\phi X}\xi - \nabla_{X}\xi  -2\phi X = 0.
\end{equation}
We may easily check that (\ref{253}) is equivalent to $(C_{2})$. and
hence
\begin{equation*}
-(\nabla_{\phi X}\eta)(\phi Y) - (\nabla_{X}\eta)(Y)-2g(\phi X,Y) =
0,
\end{equation*}
and hence
\begin{equation*}
(\nabla_{\phi X}\eta)(\phi Y) + (\nabla_{X}\eta)(Y)+2g(\phi X,Y) = 0
\end{equation*}
for any $X, Y  \in\mathfrak{X}(M)$.\; \; \;
\;\;\;\;\qquad\qquad\qquad\qquad\qquad\qquad\qquad\qquad\qquad\qquad\qquad\qquad
$\square$

% We denote by $\nabla$ the Levi-Civita connection and by $R$ the
% corresponding Riemann curvature tensor field given
% by $R(X,Y)$=[$\nabla_X$, $\nabla_Y$]- $\nabla_{[X,Y]}$ for all
% vector fields $X$,$Y$ on $M$.
% We denote by $\rho$ the Ricci tensor field of type (0,2) % by $Q$ the
%Ricci operator,
%and by $\tau$ the scalar curvature. %We define on $M$
%the operators \emph{h,~l} by
%\begin{equation}\label{23}
%        hX=\frac{1}{2}({\pounds}_\xi \phi)X,\quad lX=R(X,\xi)\xi
%\end{equation}
%where ${\pounds}_\xi$ is the Lie derivative in the direction of
%$\xi$. {{It is easily checked that {$h$ and $\l$ are symmetric
%operators}} and satisfy the following equalities
%\begin{equation}\label{24}
%        h\xi = 0, \quad {{l\xi = 0,}} \quad    h\phi = - \phi h.
%\end{equation}
%}}\\ We also have the following formulas for a contact metric
%manifold :
%\begin{equation}\label{25}
% \begin{split}
%   &\nabla_X \xi=-\phi X-\phi hX, \quad(\text{and hence }\nabla_\xi
%        \xi=0)\\
%   &\nabla_\xi \phi=0,\\
%   &  trl=g(Q\xi,\xi)=2{{n}}-tr(h^2),\\
%   &\phi l \phi-l=2(\phi^2 + h^2),\\
%   & \nabla_\xi h=\phi-\phi l-\phi h^2.\\
%\end{split}
%\end{equation}

%%%%%%%%%%%%%%%%%%%%%%%%%%%%%%%%%%%%%%%%%%%%%%%%%%%%%%%%%%%%%%%%%%%%%%%%%%%%%%%%%%%%%%%
\section{A characterization of contact metric manifolds}\label{sec3}
%%%%%%%%%%%%%%%%%%%%%%%%%%%%%%%%%%%%%%%%%%%%%%%%%%%%%%%%%%%%%%%%%%%%%%%%%%%%%%%%%%%%%%%

First of all, we shall show the following.

\begin{lem}\label{3.1}
Let $M = (M, \phi, \xi,\eta, g)$ be an almost contact metric
manifold satisfying the equality $(C_{3})$ in Proposition 2.6. Then,
the tensor field $h$ anti-commutes with $\phi$ and the following
equality\\
\begin{equation}\label{31}
g(hX,Y)-g(hY,X) = -\frac{1}{2}N^{(2)}(X,Y)
\end{equation}
holds for any $X, Y  \in\mathfrak{X}(M)$.
\end{lem}

\noindent {\bf{Proof.}} From the hypotheses ,it follows irrediately
that $M$ satisfies the equality $(C_{4})$. Thus, taking account of
(\ref{212}), we have
\begin{equation}\label{32}
\begin{split}
(\phi h + h\phi)X & = \phi h X +h \phi X \\
 & = \frac{1}{2}(-\phi\nabla_{\phi X}\xi + {\phi}^{2}\nabla_{X}\xi
 -\nabla_{{\phi}^{2}X}\xi + \phi\nabla_{\phi X}\xi)\\
& = \frac{1}{2}(-\phi\nabla_{\phi X}\xi - \nabla_{X}\xi
+ \eta(\nabla_{X}\xi)\xi + \nabla_{X}\xi - \eta(X)(\nabla_{\xi}\xi) + \phi\nabla_{\phi X}\xi) \\
& = \frac{1}{2}(\eta(\nabla_{X}\xi)\xi-\eta(X)(\nabla_{\xi}\xi))\\
& = 0
\end{split}
\end{equation}
for any $X  \in\mathfrak{X}(M)$. and hence $h$ anti-commutes with
$\phi$. Further, from (\ref{212}) with $(C_{3})$ and (\ref{210}), we
have
\begin{equation*}
\begin{split}
& g(hX,Y)-g(hY,X)\\
   % & = \frac{1}{2} \bigg( -g(\nabla_{\phi X}\xi,Y) +
   %    g(\phi\nabla_{X}\xi,Y)+g(\nabla_{\phi
   %    Y}\xi,X)-g(\phi\nabla_{Y}\xi,X) \bigg),\\
& {=\frac{1}{2} \bigg( -\nabla_{\phi X}\eta(Y) -
    (\nabla_{X}\eta)(\phi Y)+ (\nabla_{\phi
    Y}\eta)(X)+ (\nabla_{Y}\eta)(\phi X) \bigg)}\\
& = -\frac{1}{2}N^{(2)}(X,Y)
\end{split}
\end{equation*}
for any $X, Y  \in\mathfrak{X}(M)$.

Now, let $M=(M, \phi ,\xi,\eta,g)$ be a contact metric manifold.
Then, it is well-known that the tensor field $h$ is symmetric with
respect to the metric $g$ and anti-commutes with $\phi$ and $M$
satisfies the following conditions
\begin{equation*}
(C_{0}) \quad \nabla_{X}\xi = -\phi X - \phi hX,
\end{equation*}
and
\begin{equation*}
(C_{1}) \quad (\nabla_{X}\phi)Y + (\nabla_{\phi X}\phi)\phi Y =
2g(X,Y)\xi - \eta(Y)X - \eta(X)\eta(Y)\xi - \eta(Y)hX
\end{equation*}
for any $X, Y  \in\mathfrak{X}(M)$. Thus, we see that the equalities
$(C_{2}) \sim (C_{4})$ in Proposition \ref{2.6} hold on $M$ and
$(C_{0})$ is equivalent to (\ref{230})(and hence (\ref{229})) by
virtue of Proposition \ref{2.6} together with its proof.\; \; \;
\;\;\; \;
\;\;\;\;\qquad\qquad\qquad\qquad\qquad\qquad\qquad\qquad\qquad\qquad
$\square$

\noindent Thus, from the above arguments and Lemma \ref{3.1}, we
have the following theorem.

\begin{thm}\label{3.2}
A contact metric manifold is characterized as an almost contact
metric manifold $M=(M, \phi ,\xi,\eta,g)$ satisfying the following
conditions
\begin{equation*}
\begin{split}
& (C) \quad h \ is \ symmetric \\
& {(C_{1}) \quad (\nabla_{X}\phi)Y+(\nabla_{\phi X}\phi)\phi Y =
2g(XY)\xi - \eta(Y)X - \eta(X)\eta(Y)\xi - \eta(Y)hX}
\end{split}
\end{equation*}
for any $X, Y  \in\mathfrak{X}(M)$.
\end{thm}

\noindent {\bf{Proof.}} First, from Proposition \ref{2.6}, it
follows that $M$ satisfies the conditions $(C_{2}) \sim (C_{4})$ in
the some Proposition. thus, from $(C_{2})$, we have
\begin{equation*}
(\nabla_{X}\eta)(Y) + (\nabla_{\phi X} \eta)(\phi Y) = -2g(\phi
X,Y),
\end{equation*}
and hence
\begin{equation}\label{33}
(\nabla_{X}\eta)(Y) - (\nabla_{Y} \eta)( X) + (\nabla_{\phi X}
\eta)(\phi Y) -(\nabla_{\phi Y}\eta)(\phi X) = -4g(\phi X,Y)
\end{equation}
for any $X, Y  \in\mathfrak{X}(M)$. Further, from the condition
$(C)$, taking account of (\ref{210}) and (\ref{215}), we have
\begin{equation}\label{34}
(\nabla_{\phi X}\eta)(Y)-(\nabla_{Y}\eta)(\phi X) - (\nabla_{\phi
Y}\eta)(X)+(\nabla_{X}\eta)(\phi Y) = 0
\end{equation}
for any $X, Y  \in\mathfrak{X}(M)$. From (\ref{34}), we have also
\begin{equation*}
(\nabla_{{\phi}^{2}X}\eta)(Y)-(\nabla_{Y}\eta)({\phi}^{2}X)-(\nabla_{\phi
Y}\eta)(\phi X)+(\nabla_{\phi X}\eta)(\phi Y) =0,
\end{equation*}
and hence
\begin{equation}\label{35}
\begin{split}
& -(\nabla_{X}\eta)(Y) + \eta(X)(\nabla_{\xi}\eta)(Y) + (\nabla_{Y}\eta)(X)\\
& - \eta(X)(\nabla_{Y}\eta)(\xi)-(\nabla_{\phi Y}\eta)(\phi X) +
     (\nabla_{\phi X}\eta)(\phi Y) = 0
\end{split}
\end{equation}
for any $X, Y  \in\mathfrak{X}(M)$. Thus, from (\ref{33}) and
(\ref{35}), we have
\begin{equation*}
2(\nabla_{X}\eta)(Y)-2(\nabla_{Y}\eta)(X) = -4g(\phi X,Y),
\end{equation*}
and hence
\begin{equation}\label{36}
 { (\nabla_{X}\eta)(Y)-(\nabla_{Y}\eta)(X) = -2g(\phi X, Y)=
 2\Phi(X,Y),}
 \end{equation}

\noindent namely

\begin{equation*}
d\eta(X,Y) = \Phi(X,Y)
\end{equation*}
for any $X, Y  \in\mathfrak{X}(M)$. Therefore, $M$ is a contact
metric manifold. The converse is evident. This completes the proof
of Theorem \ref{3.2}. \; \; \; \;\;\;
\;\;\;\;\qquad\qquad\qquad\qquad\qquad\qquad $\square$

%%%%%%%%%%%%%%%%%%%%%%%%%%%%%%%%%%%%%%%%%%%%%%%%%%%%%%%%%%%%%%%%%%%%%%%%%%%%%%%%%%%%%%%
\section{Quasi contact metric manifolds}\label{sec4}
%%%%%%%%%%%%%%%%%%%%%%%%%%%%%%%%%%%%%%%%%%%%%%%%%%%%%%%%%%%%%%%%%%%%%%%%%%%%%%%%%%%%%%%

First of all, we shall show the following

\begin{lem}\label{4.1}
Let $M=(M, \phi, \xi, \eta , g)$ be an almost contact metric
manifold. Then the conditions $(C_{1})$ and $(C'_{1})$ are
equivalent to each other.
\end{lem}

\noindent {\bf{Proof.}} We assume that $M$ satisfies the condition
$(C_{1})$. Then, it follows from Proposition \ref{2.6} that $M$
satisfies the condition $(C_{2})$, and hence, we have
\begin{equation}\label{41}
\nabla_{X}\xi - \phi \nabla_{\phi X}\xi +2\phi X = 0
\end{equation}
for any $X, Y  \in\mathfrak{X}(M)$. Thus, from (\ref{41}) we get
\begin{equation*}
\phi \nabla_{X}\xi - {\phi}^{2}\nabla_{\phi X}\xi + 2 {\phi}^{2}X =
0,
\end{equation*}
and hence
\begin{equation}\label{42}
\phi \nabla_{X}\xi = -\nabla_{\phi X}\xi +2X -2 \eta(X)\xi.
\end{equation}

\noindent Since $M$ satisfies the condition $(C_{3})$ by virtue of
Proposition \ref{2.6}, from (\ref{212}) and (\ref{42}), we have

\begin{equation}\label{43}
\begin{split}
hX & = \frac{1}{2}(-\nabla_{\phi X}\xi + \phi \nabla_{X}\xi)\\
   & = \frac{1}{2}(-2\nabla_{\phi X}\xi +2X -2\eta(X)\xi)\\
   & = -\nabla_{\phi X}\xi + X - \eta(X)\xi.
\end{split}
\end{equation}
Thus, from (\ref{43}), we see that the equality $(C_{1})$ reduces to

\begin{equation}\label{44}
\begin{split}
& (\nabla_{X}\phi)Y + (\nabla_{\phi X} \phi)\phi Y \\
&= 2g(X,Y)\xi - \eta(Y)X -\eta(X)\eta(Y)\xi - \eta(Y)(-\nabla_{\phi
X}\xi + X -\eta(X)\xi)\\
& = 2g(X,Y)\xi -2\eta(Y)X + \eta(Y)\nabla_{\phi X}\xi
\end{split}
\end{equation}
for any $X, Y  \in\mathfrak{X}(M)$. The equality (\ref{44}) is
nothing but the equality $(C'_{1})$.

Conversely, we assume that $M$ satisfies the condition $(C'_{1})$.
Then, it follows from Proposition \ref{2.7} that $M$ also satisfies
the condition $(C_{2})$, and hence we have (\ref{41}) and hence,
(\ref{42}) and (\ref{43}). Thus, finally we see that the equality
$(C'_{1})$ reduces $(C_{1})$.  \;\;\;\; $\square$

\noindent Therefore, from Proposition 2.5 and Lemma 4.1, we have the
following Theorem.

\begin{thm}\label{4.2}
A quasi contact metric manifold is characterized as an almost
contact metric manifold $M=(M,\phi, \xi, \eta, g)$ satisfying the
following condition $(C_{1})$ :
\begin{equation*}
(C_{1}) \quad (\nabla_{X}\phi)Y+(\nabla_{\phi X}\phi)\phi Y =
2g(X,Y)\xi - \eta(Y)X-\eta(X)\eta(Y)\xi - \eta(Y)hX
\end{equation*}
for any $X, Y  \in\mathfrak{X}(M)$.
\end{thm}

\noindent { \bf{Remark} } It is well-known that a 4-dimensional
quasi K\"ahler manifold is necessarily an almost K\"ahler manifold.
Thus, a $3$-dimensional quasi contact metric manifold is necessarily
a contact metric manifold.  Some classes of 3-dimensional contact
metric manifolds have been discussed in \cite {JPS}. From our
discussion in this paper, the following question will naturally
arise.

\begin{que}\label{1}
Does there exist a $(2n+1)(n\geq2)$-dimensional quasi contact metric
manifold which is not a contact metric manifold?
\end{que}

%%%%%
% (1) $\quad$ Let
% $\bar{M}=(\bar{M},\bar{J}_{\alpha},\bar{g}_{\alpha})(\alpha=1,2,3)$
% be hyper quasi K\"ahler manifold, namely,
% $(\bar{J}_{\alpha},\bar{g}_{\alpha})(\alpha=1,2,3)$ defines a hyper
% almost Hermitian structure on $M$ and each
% $(M,\bar{J}_{\alpha},\bar{g}_{\alpha})(\alpha=1,2,3)$ is a quasi
% K\"ahler manifold. Then N.Murakirhi, k.Sekigawa-Yamada [] proved
% that $\bar{M}=(M,\bar{J}_{\alpha},\bar{g}_{\alpha})(\alpha=1,2,3)$
% is a hyper K\"ahler manifold, which is a generalization of the
% result by Hitchin {\cite{H}.

% \begin{thm}
% Any hyper almost K$\ddot{a}$hler manifold is a hyper K\"ahler
% manifold.
% \end{thm}

% Now, let
% $(M,\phi_{\alpha},\xi_{\alpha},\eta_{\alpha},g)(\alpha=1,2,3)$ be a
% $(4n-1)$-dimension contact metric manifold, namely, smooth manifold
% with contact metric 3-structure.

% \begin{thm}
% Let $M=(M,\phi_{\alpha},\xi_{\alpha},\eta_{\alpha},g)(\alpha=1,2,3)$
% be a $(4n-1)$-dimension smooth manifold equipped with a contact
% metric 3-structure
% $(\phi_{\alpha},\xi_{\alpha},\eta_{\alpha},g)(\alpha=1,2,3)$ ,Then
% \{$(\phi_{\alpha},\xi_{\alpha},\eta_{\alpha},g)\}(\alpha=1,2,3)$ is
% a Sasakian 3-structure on $M$.
% \end{thm}

%%%%%%%%%%%%%%%%%%%%%%%%%%%%%%%%%%%%%%%%%%%%%%%%%%%%%%%%%%%%%%%%%%%%%%%%%%%%%%%%%%%%%%%
\section*{Acknowledgement}
%%%%%%%%%%%%%%%%%%%%%%%%%%%%%%%%%%%%%%%%%%%%%%%%%%%%%%%%%%%%%%%%%%%%%%%%%%%%%%%%%%%%%%%
This work was supported by the National Research
Foundation of Korea (NRF) grant funded by the Korea government
(MEST) (2013020825).

%%%%%%%%%%%%%%%%%%%%%%%%%%%%%%%%%%%%%%%%%%%%%%%%%%%%%%%%%%%%%%%%%%%%%%%%%
%\section{Applications}
%%%%%%%%%%%%%%%%%%%%%%%%%%%%%%%%%%%%%%%%%%%%%%%%%%%%%%%%%%%%%%%%%%%%%%%%%

\medskip

\noindent


\begin{thebibliography}{20}

\def\refprep#1#2{#1, {\em #2,} preprint.}
\def\refnot#1#2#3{#1, #2, {\em #3},\ to appear.}
\def\refart#1#2#3#4{#1, #2, {\em #3\/} #4.}
\def\refbook#1#2#3{#1, {\em #2,\/} #3.}

%\bibitem{Ber}
%            M. Berger,
%            Quelques formules de variation pour une structure riemannienne,
%            Ann. Sci. \'Ecole Norm. Sup. 4$^{e}$ Serie 3
%            (1970),
%            285--294.

\bibitem{B1}
D. E. Blair, Riemannian geometry of contact and symplectic
manifolds. Second edition, Progress in Math. 203 (2002),
Birkh\"auser, Boston.

\bibitem{CG}  D. Chinea and C. Gonzalez,
A classification of almost contact metric manifolds, Ann. Mat. Pura
Appl. (4) 156 (1990), 15-36.




\bibitem{GH}
A. Gray, L. M. Hervella, The sixteen classes of almost Hermitian
manifolds and their linear imvarients, Ann. Mat. Pura Appl.
123 (1980), 35-58.

%\bibitem{H} N. J. Hitchin, The self-duality equations on a Riemannian surface,
%Proc. London Math. Soc. 55 (1987), 59-126.


\bibitem{JPS} J. E. Jin J. H. Park and K. Sekigawa, Notes on some classes of 3-dimensional contact metric
manifolds, Balkan J. Geom. Appl. (2012) 17(2) 42-53.

%\bibitem{K}
%T. Kashiwada, On a contact 3-structure, Math. Z. 238 (2001),
%829-832.

%\bibitem{M}
%N. Murakashi, K. Sekigawa, and A. Yamada, Integrability of almost
%quaternianie manifolds, Indian J. Math. 42 (2000), 313-329.

%\bibitem{S}
%K. Sekigawa and A. Yamada, A. Note on the integrability of a class
%of almost quaternianie manifolds, Indian J. Math. 48 (2006),
%239-248.




% \bibitem{B2}
%D. E. Blair, T. Koufogiorgos and V. J. Papantoniou, Contact metric
%manifolds satisfying a nullity condition,  Israel J. Math. 91
%(1995), 189-214.

%\bibitem{B3}
%D. E. Blair, T. Koufogiorgos and R. Sharma, A Classification of
%3-dimensional contact metric manifolds with $Q\varphi$=$\varphi Q$,
%Kodai Math. J. 13 (1990), 391-401.







%\bibitem{La}
%            M.-L. Labbi,
%            Variational properties of the Gauss-Bonnet curvatures,
%            Calc. Var.
%            32
%            (2008),
 %           175--189.

%\bibitem{T}
%Y. Tashiro, On contact structure of hypersurfaces in complex
%manifolds. I, Tohoku Math. J. 15 (1963), 62-78.

\bibitem{T2}
Y. Tashiro, On contact structure of hypersurfaces in complex
manifolds. I, II Tohoku Math. J. 15 (1963), 62-78 and 167-175.
%\bibitem{Tanno}
%S. Tanno, Variational problems on contact Riemannian manifolds,
%Trans. Amer. Math. Soc. 314 (1989), 349-379.





\end{thebibliography}
\end{document}